\documentclass[12pt]{article}
\usepackage{amssymb}
\usepackage{amsfonts}
\usepackage{latexsym}
\usepackage{amsthm}

\newcommand{\er}[1]{{\rm(\ref{#1})}}
\def\lb{\label}
\theoremstyle{plain}

\theoremstyle{remark}

\begin{document}

\def\a{\alpha}           \def\cA{{\cal A}}     \def\bA{{\bf A}}
\def\b{\beta}            \def\cB{{\cal B}}     \def\bB{{\bf B}}
\def\g{\gamma}           \def\cC{{\cal C}}     \def\bC{{\bf C}}
\def\G{\Gamma}           \def\cD{{\cal D}}     \def\bD{{\bf D}}
\def\d{\delta}           \def\cE{{\cal E}}     \def\bE{{\bf E}}
\def\D{\Delta}           \def\cF{{\cal F}}     \def\bF{{\bf F}}
\def\ve{\varepsilon}     \def\cG{{\cal G}}     \def\bG{{\bf G}}
\def\z{\zeta}            \def\cH{{\cal H}}     \def\bH{{\bf H}}
\def\e{\eta}             \def\cI{{\cal I}}     \def\bI{{\bf I}}
\def\vt{\vartheta}       \def\cJ{{\cal J}}     \def\bJ{{\bf J}}
\def\vT{\Theta}          \def\cK{{\cal K}}     \def\bK{{\bf K}}
\def\k{\kappa}           \def\cL{{\cal L}}     \def\bL{{\bf L}}
\def\l{\lambda}          \def\cM{{\cal M}}     \def\bM{{\bf M}}
\def\L{\Lambda}          \def\cN{{\cal N}}     \def\bN{{\bf N}}
\def\m{\mu}              \def\cO{{\cal O}}     \def\bO{{\bf O}}
\def\n{\nu}              \def\cP{{\cal P}}     \def\bP{{\bf P}}
\def\r{\rho}             \def\cQ{{\cal Q}}     \def\bQ{{\bf Q}}
\def\s{\sigma}           \def\cR{{\cal R}}     \def\bR{{\bf R}}
\def\S{\Sigma}           \def\cS{{\cal S}}     \def\bS{{\bf S}}
\def\t{\tau}             \def\cT{{\cal T}}     \def\bT{{\bf T}}
\def\f{\phi}             \def\cU{{\cal U}}     \def\bU{{\bf U}}
\def\F{\Phi}             \def\cV{{\cal V}}     \def\bV{{\bf V}}
\def\vp{\varphi}         \def\cW{{\cal W}}     \def\bW{{\bf W}}
\def\c{\chi}             \def\cX{{\cal X}}     \def\bX{{\bf X}}
\def\p{\psi}             \def\cY{{\cal Y}}     \def\bY{{\bf Y}}
\def\P{\Psi}             \def\cZ{{\cal Z}}     \def\bZ{{\bf Z}}
\def\o{\omega}
\def\O{\Omega}
\def\x{\xi}
\def\X{\Xi}
\def\eps{\epsilon}
\def\vk{\varkappa}

\def\Z{{\Bbb Z}}
\def\R{{\Bbb R}}
\def\C{{\Bbb C}}
\def\T{{\Bbb T}}
\def\N{{\Bbb N}}
\def\S{{\Bbb S}}
\def\H{{\Bbb H}}
\def\dD{{\Bbb D}}

\let\ge\geqslant
\let\le\leqslant
\let\geq\geqslant
\let\leq\leqslant
\def\ma{\left(\begin{array}{cc}}
\def\am{\end{array}\right)}
\def\iint{\int\!\!\!\int}
\def\lt{\biggl}
\def\rt{\biggr}
\let\geq\geqslant
\let\leq\leqslant
\def\[{\begin{equation}}
\def\]{\end{equation}}
\def\wt{\widetilde}
\def\pa{\partial}
\def\sm{\setminus}
\def\es{\emptyset}
\def\no{\noindent}
\def\ol{\overline}
\def\iy{\infty}
\def\ev{\equiv}
\def\/{\over}
\def\ts{\times}
\def\os{\oplus}
\def\ss{\subset}
\def\h{\hat}
\def\Re{\mathop{\rm Re}\nolimits}
\def\Im{\mathop{\rm Im}\nolimits}
\def\supp{\mathop{\rm supp}\nolimits}
\def\sign{\mathop{\rm sign}\nolimits}
\def\Ran{\mathop{\rm Ran}\nolimits}
\def\Ker{\mathop{\rm Ker}\nolimits}
\def\Tr{\mathop{\rm Tr}\nolimits}
\def\const{\mathop{\rm const}\nolimits}
\def\dist{\mathop{\rm dist}\nolimits}
\def\diag{\mathop{\rm diag}\nolimits}
\def\Wr{\mathop{\rm Wr}\nolimits}
\def\BBox{\hspace{1mm}\vrule height6pt width5.5pt depth0pt \hspace{6pt}}

\def\Twelve{
\font\Tenmsa=msam10 scaled 1200
\font\Sevenmsa=msam7 scaled 1200
\font\Fivemsa=msam5 scaled 1200
\textfont\msbfam=\Tenmsb
\scriptfont\msbfam=\Sevenmsb
\scriptscriptfont\msbfam=\Fivemsb

\font\Teneufm=eufm10 scaled 1200
\font\Seveneufm=eufm7 scaled 1200
\font\Fiveeufm=eufm5 scaled 1200
\textfont\eufmfam=\Teneufm
\scriptfont\eufmfam=\Seveneufm
\scriptscriptfont\eufmfam=\Fiveeufm}

\def\Ten{
\textfont\msafam=\tenmsa
\scriptfont\msafam=\sevenmsa
\scriptscriptfont\msafam=\fivemsa

\textfont\msbfam=\tenmsb
\scriptfont\msbfam=\sevenmsb
\scriptscriptfont\msbfam=\fivemsb

\textfont\eufmfam=\teneufm
\scriptfont\eufmfam=\seveneufm
\scriptscriptfont\eufmfam=\fiveeufm}

\title {The Lyapunov function for Schr\"odinger operators with a periodic  $2\ts 2$ matrix potential  }

\author{Andrei Badanin
\begin{footnote}
{ Department of
 Mathematics of Archangel University, Russia e-mail: badanin@agtu.ru}
\end{footnote}
\and Jochen Br\"uning
\begin{footnote} { Institut f\"ur
Mathematik, Humboldt Universit\"at zu Berlin, e-mail:
bruening@math.hu-berlin.de}
\end{footnote}
 \and Evgeny Korotyaev
\begin{footnote}
{corresponding author, Institut f\"ur  Mathematik,  Humboldt Universit\"at zu Berlin, Rudower Chaussee 25, 12489, Berlin, Germany
e-mail:
 evgeny@math.hu-berlin.de\ \  }
\end{footnote}
}

\maketitle

\begin{abstract}
\no We consider the Schr\"odinger operator on the real line with a
$2\ts 2$ matrix valued  1-periodic potential. The spectrum of this
operator is absolutely continuous and consists of intervals
 separated by gaps. We define a Lyapunov function
which is analytic on a two sheeted Riemann surface. On each
sheet, the Lyapunov function has the same properties 
as in the scalar case, but it has 
branch points, which we call resonances. We prove the existence of real as well as non-real resonances for specific potentials. We determine the asymptotics of the periodic and anti-periodic spectrum and of the resonances at high energy. We show that there exist two type of
gaps: 1) stable gaps, where the endpoints are periodic and
anti-periodic eigenvalues, 2) unstable (resonance) gaps, where the
endpoints are resonances (i.e., real branch points of the Lyapunov
function). We also show that
periodic and anti-periodic spectrum together determine the spectrum of the matrix Hill operator.
\end{abstract}


\section {Introduction and main results}
\setcounter{equation}{0}

We consider the self-adjoint operator $Ty=-y''+V(x)y,$ acting in
$L^2(\R)\os L^2(\R),$ where $V$ is a symmetric 1-periodic $2\ts 2$ matrix potential which belongs to the real space  $\cH^p, p=1,2$, given by
$$
\cH^p=\lt\{V=V^*=V^T=\ma V_1&V_3\\
V_3&V_2\am :\ \ \ \ \ \ \int_0^1V_3(x)dx=0 \rt\},\ \ \
$$
equipped with the norm $\|V\|_p^p=\int_0^1
(|V_1(x)|^p+|V_2(x)|^p+2|V_3(x)|^p)dx<\iy$.
Without loss of generality {\bf we assume} 
$$
V_{(0)} = \int_0^1V(t)dt = {\rm diag} \{V_{10}, V_{20}\},\ \ \ \ V_{10}\le V_{20},\ \  V_{m0}=\int_0^1V_{m}(x)dx,\  \ m=1,2.
$$
Introduce the self-adjoint operator
$T^0=-{d^2\/dx^2}$, with the domain $Dom(T^0)=W_2^2(\R)\os W_2^2(\R)$. In order to get self-adjointness of $T$ we use symmetric quadratic forms. We shortly recall a well known argument (see [RS]). We define the quadratic form
$
(V\p, \p_1)=-\int_\R V\p\ \ol\p_1 dx,
\ \ \ \ \ \p,\p_1\in Dom (T^0).
$
Using the estimate (see [K1]) $|(q'f,f)|< \ve (f',f)+b_\ve(f,f)$ for any small $\ve >0$ and some $b_\ve>0$ and any $f\in W_2^2(\R),q\in L^2(\R/\Z)$ we deduce that
$$
|(V\p, \p_1)|<(1/2)(\p', \p')+b(\p, \p), \ \ \ \p\in W_2^2(\R)\os W_2^2(\R).
$$
Thus we can apply the KLMN Theorem (see [RS]) to define
$T=-{d^2\/dx^2}+V$. There exists a unique self-adjoint operator $T$ with form domain $\cQ (T)=W_1^2(\R)\os W_1^2(\R)$ and
$$
(T\p, \p_1)=(-\p'', \p_1)+(V\p,\p_1),\ \ \ \ {\rm all}\ \ \ \
\p, \p_1\in \cQ (T^0)=W_1^2(\R).
$$
Any domain of essential self-adjointness for $T^0$ is a form core for $T$.

It is well known (see [DS] p.1486-1494, [Ge]) that the spectrum
$\s(T)$ of $T$ is absolutely continuous and consists of
non-degenerated intervals $S_n,n=1,2..,$. These intervals are separated by the gaps $G_n$ with lengths $|G_n|> 0, n=1,2..,N_G\le \iy$.  Introduce the fundamental $2\ts 2$-matrix solutions $\vp(x,\l)$, $\vt(x,\l)$ of the equation
\[
 \lb{1b}
 -y''+V(x)y=\l y,\ \ \ \l\in\C,
\]
with initial conditions $\vp(0,\l)=\vt'(0,\l)=0$,
$\vp'(0,\l)=\vt(0,\l)=I_2$, where $I_m,m\ge 1$ is the identity
$m\ts m$ matrix.  Here and below we use the notation
$(')=\pa /\pa x$. We define the $4\ts 4$-monodromy matrix by
\[
\lb{1c}
M(\l)=\cM(1,\l), , \ \ \ \ \ \ 
\cM(x,\l)=\ma \vt(x,\l)&\vp(x,\l)\\
\vt'(x,\l)&\vp'(x,\l)\\ \am.
\]
The matrix valued function $M$ is entire. An eigenvalue of $M(\l)$
is called a {\it multiplier}, i.e., it is a root of the algebraic
equation $D(\t,\l)=0$, where 
\[
\lb{1det}
D(\t,\l):=\det (M(\l)-\t I_4), \ \ \t,\l\in\C.
\]
There is an enormous literature on the scalar Hill operator
including the inverse spectral theory [M], [GT], [KK]. In the
recent paper [K1] one of the authors solved the inverse problem
(including characterization) for the operator $-y''+v'y$, where a
function $v\in L^2(\T)$ and $v'$ is the distribution and
$\T=\R/\Z$. However, in spite of the importance of extending
these studies to vector differential equations,
apart from the information given
by Lyapunov and Krein (see [YS]), until recently nothing essential  has been done. The matrix potential poses interesting new problems: 1) to construct the Lyapunov function, 2) to define the quasimomentum as a conformal mapping, 3) to derive appropriate trace formulae (e.g. analog to the trace formulas in [K2]), 4) to obtain a priori estimates of potentials in terms of gap lengths, 5) to define and to study the integrated density of states. In fact this is the motivation of our paper.

 The basic results in the direct spectral theory for the
matrix case were obtained by Lyapunov [Ly] and Krein [Kr] (see
also  Gel'fand and Lidskii [GL]).
 Below we need the following well-known results of Lyapunov [YS],
 which we formulate only for the case $2\ts 2$ matrices.

{\bf Theorem (Lyapunov )} {\it Let $V\in \cH^1$. 

\no  i) If $\t(\l)$  is a multiplier for some $\l\in\C$
(or $\l\in\R$), then
$\t^{-1}(\l)$ (or $\ol\t(\l)$ ) is a multiplier too.

\no ii) $M(\l)$ has exactly four multipliers
$\t_1(\l),\t_1^{-1}(\l),\t_2(\l),\t_2^{-1}(\l)$ for all $\l\in\C$. Moreover, $\l\in\s(T)$ iff $|\t_m(\l)|=1$ for some $m=1,2$.

\no iii) If $\t(\l)$ is a simple multiplier and $|\t(\l)|=1$, then
$\t'(\l)\ne 0$}.

We mention the papers relevant in our context,
for various problems periodic systems we refer the reader to [CGHL].
In the paper [Ca1,2] Carlson  obtained the trace formulas.
 In [Ca3] he proved the compactness of
Floquet isospectral sets for the matrix Hill operator.

By the Lyapunov Theorem, each $M(\l), \l\in\C$ has exactly four
multipliers $\t_m(\l),\t_m^{-1}(\l)$, $m=1,2$, which are the roots
of the characteristic polynomial $\det (M(\l)-\t I_4)=0$. If
$c_0={V_{20}-V_{10}\/2}>0$, then the multipliers have the following asymptotics
\[
 \lb{1ca}
\t_m(\l)=e^{i(\sqrt \l-{V_{m0}\/2\sqrt \l}+O(1/\l))},\ \ \
|\l|\to\iy,\ \ \  |\sqrt \l-\pi n|>{\pi\/4},\ m=1,2,
\]
(see  Lemma 3.2). Next we define the functions
\[
\lb{1d}\lb{dv} 
 \m_m(\l)={\Tr\/4} M^m(\l),\ \ \ m=1,2,\ \ \ \
 \ \
  \r={\m_2+1\/2}-\m_1^2,\ \ \ \
\r_0(\l)=c_0{\sin\sqrt \l\/2\sqrt \l}.
\]
Note that $\vp(x,\l), \vt(x,\l), \m_1(\l),\m_2(\l),\r(\l)$ are real for all  $x,\l\in \R$ and entire. If
$c_0={V_{20}-V_{10}\/2}>0$, then Lemma 3.1 yields $\r(\l)=\r_0^2(\l)(1+o(1))$, as $|\l|\to \iy$ in $\cD_1$, where
\[
\lb{DD} \cD_r:=\rt\{\l\in \C: |\l|>r^2, |\sqrt \l-\pi
n|>{\pi\/4}, n\in \N \rt\},\ \ r>0.
\]
 Thus we may define the analytic function
$\sqrt{\r(\l)}, \l\in \cD_r$, for large  $r>0$, by the requirement 
that $\sqrt{\r(t)}=\r_0(\l)(1+o(1))$ as $|\l|\to\iy$ in $\cD_r$.
 Then there exists a unique analytic continuation of $\sqrt
\r$ from $\cD_r$ into the two sheeted Riemann surface $\L$
(in general, of infinite genus) defined by $\sqrt \r$. We now introduce our Lyapunov function $\D(\l)$ by
\[
\lb{1g}
 \D(\l)=\m_1(\l)+\sqrt{\r(\l)},\ \ \ \ \ \ \l\in\L.
\]
Let $\D_1(\l)=\m_1(\l)+\sqrt{\r(\l)}$ on the first sheet $\L_1$
and let $\D_2(\l)=\m_1(\l)-\sqrt{\r(\l)}$ on the second sheet
$\L_2$ of $\L$.
Now we formulate our main result concerning the function $\D$.

\no {\bf Theorem 1.1.} {\it Let $V\in \cH^1$ with $V_{10}<V_{20}$. Then the function $\D=\m_1+\sqrt{\r}$ is analytic on the two sheeted
Riemann surface $\L$ and has the following properties:

\no i) 
\[
\D_m(\l)={\t_m(\l)+\t_m^{-1}(\l)\/2}, \ \ \ \l\in\L_m,\ \ \ m=1,2,
\lb{1ga}
\]
\[
\lb{asL}
 \D_m(\l)=\cos \sqrt\l+V_{m0}{\sin
\sqrt\l\/2\sqrt\l}+O\lt({e^{|\Im \sqrt \l|}\/\l}\rt), \ \  \
m=1,2,\  \l\in\cD_1.
\]
\no ii) $\l\in \C$ belongs to $\s(T)$ iff
$\D_m(\l)\in [-1,1]$ for some $m=1,2$.

\no iii) If $\l\in\s (T)$, then $\r(\l)\ge 0$.

\no iv) (The monotonicity property.) Let $\D_m$ be real analytic on some interval $I=(\a_1,\a_2)\ss\R$ and $-1<\D_m(\l)<1$, for any $\l\in I$ for some $m\in \{1,2\}$. Then $\D_m'(\l)\ne 0$ for each $\l\in I$ .}

\no {\bf Remark.}
i) For the scalar Hill operator the monodromy matrix
has exactly two eigenvalues $\t$ and $\t^{-1}$. The Lyapunov function ${1\/2}(\t+\t^{-1})$  is an entire function of the spectral
parameter and it defines the band-gap structure of the spectrum. By
Theorem 1.1, our Lyapunov function for the matrix Hill
operator also defines the band-gap structure of the spectrum, but
it is the sheeted analytic function.

\no ii) Consider the case of a diagonal potential, i.e. $V_3=0$. 
Then the Riemann surface degenerates into two components, 
and we get
\[
\m_1={1\/2}(\D_{(1)}+\D_{(2)}),\ \ \
\r={1\/4}(\D_{(1)}-\D_{(2)})^2,\ \ \ \ \
\sqrt{\r}={1\/2}(\D_{(1)}-\D_{(2)}),
\]
where $\D_{(m)}$ is the Lyapunov function for the scalar Hill
operator $-y''+V_my,\ m=1,2$. Thus 
$\D_1=\m_1+\sqrt{\r}=\D_{(1)}$ and $\D_2=\m_1-\sqrt{\r}=\D_{(2)}$.

iii) A lot of papers are devoted to the resonances for the Schr\"odinger operator with
compactly supported potentials on the real line, see [BKW],[K3],[S],[Z]. Assume that we have the couple constant before
the potential. If  this constant changes
then roughly speaking some resonances create the eigenvalues.
In our case, if the couple 
constant (before the periodic matrix potential) changes,
then roughly speaking some resonances create the gaps,
see Proposition 1.3.

 Let $D_{\pm}(\l)={1\/4}D(\pm 1,\l)$. The zeros
of $D_{+}(\l)$ and $D_{-}(\l)$ are the eigenvalues of  periodic
and anti-periodic problem associated to the equation $-y''+Vy=\l y$. Denote by
$\l_{2n,k}, n=0,1,..$, and $ k\in \{1,2,3,4\}$ the sequence of zeros of
$D_+$ (counted with multiplicity) such that $\l_{0,1}\le
\l_{0,2}\le \l_{2,1}\le \l_{2,2}\le\l_{2,3}\le \l_{2,4} \le
\l_{4,1}\le \l_{4,2}\le\l_{4,3}\le \l_{4,4}\dots $. Denote by
$\l_{2n-1,k}, (n,k)\in\N\ts\{1,2,3,4\}$ the sequence of zeros of
$D_-$ (counted with multiplicity) such that $\l_{1,1}\le
\l_{1,2}\le\l_{1,3}\le \l_{1,4} \le \l_{3,1}\le
\l_{3,2}\le\l_{3,3}\le \l_{3,4}\dots $. Note that
$\l_{n,k},n=0,1,..$, and $ k\in \{1,2,3,4\}$ are the eigenvalues of problem with period 2 for the equation $-y''+Vy=\l y$.

 Denote by $\{r_{n}\}_{1}^{\iy }$ the sequence of zeros 
 of $\r$ in $\C$ (counted with multiplicity) such that $0\leq |r_{1}|\le |r_{2}|\le
|r_3|\le  \dots $. We call these zeros of $\r$ the resonances of
$T$. We formulate the theorem about the recovering the spectrum of
$T$ and the asymptotics of the periodic and anti-periodic
eigenvalues and resonances at high energy. Furthermore, we write
$a_n=b_n+\ell^2(n)$ iff the sequence $\{a_n-b_n\}_{n\ge
1}\in\ell^2$. Recall that $V_0,V_{m0},m=1,2$ are defined by
\er{dv}.

\no {\bf Theorem 1.2.} {\it  Let $V\in \cH^2$ with $V_{20}-V_{10}>0$. Then the following asymptotics are fulfilled:
 \[
 \lb{a8}
\l_{n,m+k}=(\pi n)^2+V_{m0}+\ell^2(n),\ \ \ m,k=1,2,\ \ \ n\to\iy,
\]
 \[
 \lb{a9}
r_{2n-m}=(\pi n)^2+\ell^2(n),\ \ m=0,1, \ \ \ n\to\iy.
\]
\no ii) Let $V\in \cH^2, V_{20}-V_{10}>0$. Then the following statements
hold

\no a) The periodic spectrum and the anti-periodic spectrum recover
the resonances and the spectrum of the operator $T$.

\no b) The periodic (anti-periodic) spectrum is recovered by the
anti-periodic (periodic) spectrum and resonances. }

Less precise asymptotics for the case $V\in C^2$ were obtained
in [Ca1].

\no{\bf Example.} Consider the operator
$T_{\g,\n}=-{d^2\/dx^2}+q_{\g,\n}$,  where $q_{\g,\n}=aJ+\g v_\n(x)J_1$ is a potential for some constants $a,\g\in \R$, $\ J=\ma 1&0\\0&-1\\ \am, J_1=\ma 0&1\\1&0\\ \am,$
and $q_{\g,\n}$   satisfies

\no {\bf Condition A.} {\it ${a\/2\pi^2}-n_a\in (0,1)$ for some
integer $n_a\ge 0$. Each function $v_\n\in L^1(\T),\n=1, {1\/2},
{1\/3},..,$ is such that $v_\n(x)=v_\n(1-x),\ x\in (0,1)$, and
$\int_0^1v_\n(x)dx=1$ and for any $f\in C(0,k), k\in \N$, the
following convergence holds: }
\[
\lb{5a}
 \int_0^kv_\n(x)f(x)dx\to \int_0^k \d_{per}(x)f(x)dx\ \ \ {\rm as} \ \ \n\to 0,\ \ \
\d_{per}=\sum_{-\iy}^\iy\d(x-n-{1\/2}).
\]

 If $\g=0$, then we have the operator
$T^0=-y''+q^0$ with a constant potential $q^0=aJ$. In this case
there are no gaps in the spectrum. The fundamental solutions of
Eq.\er{1b} with $q^0=aJ$ have the form $\vt^0=\diag(c_+,c_-)$,
$\vp^0=\diag(s_+,s_-)$, where
$$
c_{\pm}(x,\l)=\cos \e_{\pm}x,\ \ \ s_{\pm}(x,\l)={\sin
\e_{\pm}x\/\e_{\pm}},\ \ \ \e_{\pm}=\sqrt{\l\mp a},\ \ \
\diag(b,c)=\ma b&0\\0&c\\\am
$$
and the branch of $\sqrt {\l}$ is given by $\sqrt 1=1$. Below we
sometimes will write $\r(\l,V), M(\l,V),..$, instead of $\r(\l),
M(\l) , ..$, when several potentials are being dealt with. Then
the functions $\r^0(\cdot)=\r(\cdot,q^0), \m_m^0(\cdot)=\m_m(\cdot,q^0),..$
corresponding to $q^0$ have the forms
\[
\lb{m0} \m_m^0(\l)={c_+(m,\l)+c_-(m,\l)\/2},\ \ m=1,2, \ \ \ \
\r^0(\l)={(c_+(1,\l)-c_-(1,\l))^2\/4},
\]
\[
\lb{L0} D_{\pm}^0(\l)=(1\mp c_+(1,\l))(1\mp c_-(1,\l)), \ \ \
\D_{(1)}^0(\l)=c_+(1,\l),\ \ \ \ \D_{(2)}^0(\l)=c_-(1,\l).
\]
Here $\D_{(m)}^0$ is the Lyapunov function for the Eq.
$-y''-(-1)^m ay=\l y, m=1,2$.  The entire function $c_+(1,\l)-c_-(1,\l)$ has only the simple
roots $z_n^0$ given by
\[
\lb{4c} z_n^0=(\pi n)^2+{a^2\/(2\pi n)^2}=\mp \ a+(\pi
n\pm{a\/2\pi n})^2,\ \ \  n\ge 1.
\]
Note that $z_1^0>z_2^0>...>z_{n_a}^0$ and
$z_{n_a+1}^0<z_{n_a+2}^0<...$. Thus all roots of
$\r^0=\r(\cdot,q^0)$ have multiplicity 2 and are given by \er{4c}.
 If ${2a\/\pi^2}\notin \N$,
then the zeros of $\D_{(m)}^0(\l)^2-1$ have multiplicity 2 and
are given by
\[
\l_{m,n}^{0}=a-(-1)^m(\pi n)^2,\ \ \  \  \ \l_{1,n}^{0}\neq
\l_{2,s}^{0}, \ \ m=1,2,\ n,s\ge 0.
\]
Note that $\l_{m,n}^{0}, m=1,2, n\ge 0$ are the roots of $D_\pm^0
(\l)$. If $\g, \n\neq 0$ are small enough, then there exist gaps
in the spectrum of $T_{\g,\n}$.  Define the disk
$\dD_r=\{\l:|\l|<\pi ^2 r^2\},r>0$. We show that there exist the
non-degenerated resonance gaps for some $V$. In this example some
resonances $r_n\in \R$ and some $r_n\notin \R$.

\no {\bf Proposition 1.3.} {\it Let a potential $q_{\g,\n}=aJ+\g
v_\n(x)J_1, a>0, \g\in \R$, satisfy Condition A,
and let ${4a\/\pi^2}\notin \N$. Then for any $N\ge 1+a$ there
exist small $\ve, \ve_1, \ve_2>0$ such that for any $(\g,\n)\in
(0,\ve_1)\ts(0,\ve_2)$ all zeros $z_n^\pm(q_{\g,\n})$ of
$\r(\l,q_{\g,\n})$ and the zeros $\l_{m,n}^\pm(q_{\g,\n})$ of
$\D_m^2(\l,q_{\g,\n})-1$ in the disk $\dD_{N+{1\/2}}$ are simple
and have the properties
\[
z_n^\pm(V_{\g,\n})\in \C_\pm, 1\le n\le n_a, \ \ \ \ \
z_n^-(q_{\g,\n})<z_n^+(q_{\g,\n}),\ \ n_a<n\le N,\ \ \
|z_n^\pm(q_{\g,\n})-z_n^0|<\ve_2,
\]
\[
\l_{m,n}^-(q_{\g,\n})<\l_{m,n}^+(q_{\g,\n}),\ \ \ \
|\l_{m,n}^\pm(q_{\g,\n})-\l_{m,n}^0|<\ve_2,\ \ \ n=0,..,N, m=1,2.
\]
There are no other roots of $\r(\l,q_{\g,\n})$ and
$\D_m^2(\l,q_{\g,\n})-1$ in the disk $\dD_{N+{1\/2}}$.

\no Remark.} i) If $0<a<2\pi^2$, then $\r(\l,q_{\g,\n})$ has only
real roots $z_1^\pm(q_{\g,\n})$ in the disk $\dD_{N+{1\/2}}$ for
small $\g, \n$. If $a>2\pi^2$, then $\r(\l,q_{\g,\n})$ has at
least two non-real roots $z_1^\pm(q_{\g,\n})$ for small $\g, \n$.

ii) We show that the operator $T_{\g,\n}$ has the gaps associated with the periodic or anti-periodic spectrum. Moreover, we show the
existence of new gaps (so-called resonance gaps). The endpoints of
the resonance gap are the branch points of the Lyapunov function,
and, in general, they are not the periodic (or anti-periodic)
eigenvalues. These endpoints are not stable. If they are real
($n=n_a+1,.., N$), then we have a gap. If they are complex ($1\le
n\le n_a$), then we have not a gap, we have only the branch points
of the Lyapunov function in the complex plane.

\section {Fundamental solutions}
\setcounter{equation}{0}

In this section we study $\vt,\vp$. We begin with some notational
convention. A vector $h=\{h_n\}_1^N\in \C^N$ has the Euclidean
norm $|h|^2=\sum_1^N|h_n|^2$, while a $N\ts N$ matrix $A$ has the
operator norm given by $|A|=\sup_{|h|=1} |Ah|$. The function
$\vp$ satisfies the following integral equations
\[
 \lb{2ic}
\vp(x,\l)= \vp_0(x,\l)+\int_0^x{\sin
\sqrt\l(x-t)\/\sqrt\l}V(t)\vp(t,\l)dt,\ \ \ \vp_0(x,\l)={\sin
\sqrt\l x\/\sqrt\l}I_2,\ \ \
\]
where $(x,\l)\in\R\ts\C$. The standard iterations in \er{2ic}
yields
\[
\lb{2id}
 \vp(x,\l)={\sum}_{n\ge 0} \vp_n(x,\l)\,, \quad
\vp_{n+1}(x,\l)= \int_0^x{\sin
\sqrt\l(x-t)\/\sqrt\l}V(t)\vp_n(t,\l)dt.
\]
The similar expansion $\vt={\sum}_{n\ge 0} \vt_n$ with
$\vt_0(x,\l)=(\cos \sqrt\l x) I_2$ holds. We introduce the functions
\[
\lb{di} I_m^0(\l)=\int_0^mdt\int_0^t\cos
\sqrt\l(m-2t+2s)F(t,s)ds,\ \ \ F(t,s)=\Tr V(t)V(s),\ \ \
\]
$ m=1,2.$ In Lemma 2.1 we shall show the simple identity
\[
\lb{iI} \lb{di1}
I_2^0(\l)=4I_1^0(\l)\cos \sqrt \l .
\]
We define
$|\l|_1\ev\max\{1,|\l|\}$ and
\[
\lb{dvm} V_{(0)}=\int_0^1V(x)dx,\ \ V_{(1)}=\Tr V_{(0)},\ \ \
V_{(2)}=\Tr V_{(0)}^2,\ \ \   A= e^{|\Im\sqrt\l|+\vk},\ \ \
\vk={\|V\|_1\/\sqrt{|\l|_1}}.
\]
We prove

\no {\bf Lemma 2.1.} {\it For each $(x,V)\in \R_+\ts\cH^1$ the
functions $\vp(x,\cdot),\vt(x,\cdot)$ are entire and for any $N\ge
-1$ the following estimates are fulfilled:
$$
\max\lt\{\lt|\vt(x,\l)-\sum_0^N\vt_n(x,\l)\rt|,
\lt|\sqrt\l\lt(\vp(x,\l)-\sum_0^N\vp_n(x,\l)\rt)\rt|,
\lt|{1\/\sqrt\l}\lt(\vt'(x,\l)-\sum_0^N\vt_n'(x,\l)\rt)\rt|,
$$
\[
\lt|\vp'(x,\l)-\sum_0^N\vp_n'(x,\l)\rt|\rt\} \le
{(x\vk)^{N+1}\/(N+1)!}A^{x}\ . \lb{2if}
\]
Moreover, each $\m_m(\l), m=1,2$ is real for $\l\in \R$ and entire and the following estimates are fulfilled:
\[
\lb{em01}
 |\m_m(\l)|\le A^m,\ \ \ \ \ \
 |\m_m(\l)-\cos mz|\le {m\vk}A^m,
\]
\[
\lb{em2}
 |\m_m(\l)-\cos mz-{\sin mz\/4z}mV_{(1)}|\le{(m\vk)^2\/2}A^m,
\]
\[
\lb{em}
 |\m_m(\l)-\cos mz-{\sin mz\/4z}mV_{(1)}-{1\/8z^2}
(I_m^0(\l)-{m^2\cos mz \/2}V_{(2)})|\le {(m\vk)^3\/3!}A^m,
\]
where $I_1^0,I_2^0$ satisfy \er{iI} and $z=\sqrt \l$.

 \no Proof.} We prove the estimates
of $\vp$, the proof for $\vp',\vt,\vt'$ is similar. \er{2id} gives
\[
\vp_n(x,\l)= \int\limits_{0< x_1< x_2< ...<
x_{n+1}=x}\!\!\!\!\vp_0(x_1,\l) \lt(\prod\limits_{1\le k\le
n}^{\curvearrowleft}
{\sin(\sqrt\l(x_{k+1}-x_k))\/\sqrt\l}V(x_k)\rt)dx_1dx_2...dx_n,
\lb{2ig}
\]
where for matrices $a_1,a_2,...,a_n$ we denote $\prod\limits_{1\le
k\le n}^{\curvearrowleft}a_k=a_na_{n-1}...a_1$. Substituting the
estimate $|\sqrt\l\vp_0(x,\l)|\le e^{|\Im \sqrt\l|x}$ into
\er{2ig} we obtain $| \sqrt\l \vp_n(x,\l)|\le
{(x\vk)^n\/n!}e^{|\Im z|x}, $ which shows that for each $x\ge 0$
the series \er{2id} converges uniformly on bounded subsets of $\C$.
Each term of this series is an entire function. Hence the sum is
an entire function. Summing the majorants we obtain estimates
\er{2if}.

\no We have $4\m_m=\Tr M^m(\l)=\Tr M(m,\l)=\Tr \sum_{n\ge
0}M_n(m,\l)$, where $m=1,2$ and
\[
\lb{itn} \Tr M_0(m,\l)=4\cos mz,\ \ \ \ \Tr M_n(m,\l)=\Tr
\vt_n(m,\l)+\Tr \vp_n'(m,\l),\ \  n\ge 1.
\]
The estimates $|\vp_n'(m,\l)|\le {(m\vk)^n\/n!}e^{|\Im z|m}$ and
$|\vt_n(m,\l)|\le {(m\vk)^n\/n!}e^{|\Im z|m}$ yield
\[
|\Tr M_n(m,\l)|\le 4{(m\vk)^n\/n!}e^{m|\Im\sqrt\l|},\ \ n\ge 0.
\]
Using \er{itn} we obtain
\[
\Tr M_1(m,\l)={1\/z}\int_0^m\!\!\! (\sin z(m-t)\cos zt+\cos
z(m-t)\sin zt )\Tr V(t)dt={\sin mz\/z}mV_{(1)},
\]
 and
$$
\Tr M_2(m,\l)={1\/z^2}\int_0^m\int_0^t\!\!\! \sin
z(t-s)z(m-t+s)F(t,s)dtds
$$
$$
={1\/2z^2}\int_0^m\int_0^t\!\!\! (\cos z(m-2t+2s)-\cos
zm)F(t,s)dtds= {1\/2z^2} (I_m^0(\l)-{m^2\/2}\cos mz V_{(2)}),
$$
since
$$
\int_0^m\int_0^t\!\!\! F(t,s)dtds={1\/2}\Tr
\lt(\int_0^mV(t)dt\rt)^2= {m^2\/2}V_{(2)}.
$$
  We have that $\m_1,\m_2$ are entire.
Moreover, the trace of the monodromy matrix is a sum of its
eigenvalues. By Lyapunov Theorem (see Sect.1), the set of these
eigenvalues is symmetric with respect to the real axis, as
$\l\in\R$. Hence $\m_1,\m_2$ are real on $\R$.

We show \er{di1}.
 Let $I_m^0=I_m^0(\l)$. We have $2\cos z I_1^0=Y_0+Y_1$, where
$$
Y_0=\int_0^1\int_0^t\cos 2z(t-s) F(t,s)dtds,\ \ \ \
Y_1=\int_0^1\int_0^t\cos 2z(1-t+s) F(t,s)dtds.
$$
We get $I_2^0=Y_1+Y_2+Y_3$, where
$$
Y_2=\int_1^2\int_0^1\cos 2z(1-t+s) F(t,s)dtds,\ \ \
Y_3=\int_1^2\int_1^t\cos 2z(1-t+s) F(t,s)dtds
$$
and using the new variable $\t=t-1$ we get $
Y_2=\int_0^1\int_0^\t\cos 2z(\t-s) F(\t,s)d\t ds=2Y_0. $ We use
the new variables $\t=t-1, \s=s-1$ and obtain $Y_3=
\int_0^1\int_0^\t\cos 2z(1-\t+\s) F(\t,\s)d\t d\s=Y_1, $ which
yields $I_2^0=2Y_1+2Y_0$. Thus we have \er{di1}. $\BBox$

We need the basic properties of the monodromy matrix.

\no {\bf Lemma 2.2} {\it  Let $V\in \cH^1$. Then the function
 $D(\t,\l)=\det(M(\l)-\t
I_4),\l,\t\in \C^2$ is entire on $\C^2$ and the following
identities are fulfilled:
\[
\lb{2j}
 {D_\t'(\t,\l)}=-D(\t,\l)\Tr\lt(M(\l)-\t I_4\rt)^{-1},
\]
$$
D(\t,\cdot)=\t^4-4\m_1\t^3+2(4\m_1^2-\m_2)\t^2-4\m_1\t+1
$$
\[
\lb{2l}
 =\lt(\t^2-2(\m_1-\sqrt{\r})\t+1\rt)
\lt(\t^2-2(\m_1+\sqrt{\r})\t+1\rt).
\]
\no Proof.} Let $D(\t)\equiv D(\t,\l)$. The standard identity
${D'(\t)}=D(\t)\Tr((M-\t)^{-1}{d(M-\t)\/d\t})$ yields \er{2j}. We
prove \er{2l}. Since $\det M=D(0)=1$, we have: $D(\t)=1+a \t+b
\t^2+c \t^3+d \t^4, a,b,c,d\in\C$.  Then
$D(\t)=(\t-\t_1)(\t-\t_1^{-1})(\t-\t_2)(\t-\t_2^{-1})$, where
$\t_1,\t_1^{-1},\t_2,\t_2^{-1}$ are the multipliers. Therefore
$d=1,a=c$. Then we have
\[
\lb{2u}
D(\t)=1+D'(0)\t+{1\/2}D''(0)\t^2+D'(0)\t^3+\t^4.
\]
Identity \er{2j} yields $ D'(0)=-\Tr M^{-1}=-\Tr M=-4\m_1$.
Differentiating \er{2j} we obtain $
D''(0)=-D'(\t)\Tr(M-\t)^{-1}-D(\t)\Tr(M-\t)^{-2}|_{\t=0}=4(4\m_1^2-\m_2)$.
Substituting these identities into \er{2u} we obtain the first
identity in \er{2l}. The second identity is proved by direct
calculation. $\BBox$

\section {The Lyapunov function}
\setcounter{equation}{0}

We need some results about the functions $\r,\
\r_0=c_0{\sin \sqrt\l\/2\sqrt\l}$, where $c_0={1\/2}(V_{20}-V_{10})$.

\no {\bf Lemma 3.1.} {\it i) For each $V\in \cH^1$ the function
$\r={1\/2}(\m_2+1)-\m_1^2$ is entire and real on the real line. Moreover, the following estimate is fulfilled:
\[
\lb{3a}
 |\r(\l)-\r_0^2(\l)|\le 2\vk^3e^{2|\Im\sqrt\l|+2\vk},\ \ \
\l\in\C,\ \ \ \vk={\|V\|_1/\sqrt{|\l|_1}}.
\]
\no ii) Let $V\in \cH^1$ with $c_0={1\/2}(V_{20}-V_{10})>0$. Then for each integer
$N>{2^7\|V\|_1^3\/c_0^2}$ the function $\r(\l)$ has exactly $2N$
roots, counted with multiplicity, in the disk
$\{\l:|\l|<\pi^2(N+{1\/2})^2\}$ and for each $n>N$, exactly two
roots, counted with multiplicity, in the domain $\{\l:|\sqrt\l-\pi
n|<1\}$. There are no other roots.

\no iii) Let $V\in \cH^1,c_0>0$. Then the function $\sqrt{\r}$ is
an analytic function in the domain $\cD_r, r={2^9\|V\|_1^3\/c_0^2}$ given by \er{DD} and the following estimate is fulfilled
\[
\lb{3ad}
 |\sqrt{\r(\l)}-\r_0(\l)|\le
{3C_0\/5}{|\r_0(\l)|\/\sqrt{|\l|}},\ \ \ C_0\ev
4^4{\|V\|_1^3\/c_0^2}<{\sqrt{|\l|}\/2},\ \  \l\in \cD_r.
\]
\no Proof.} i) By Lemma 2.2, $\r$ is entire and real on the real line. Let $\m_m\ev \m_{m0}+\m_{m1}+\m_{m2}+\wt\m_{m3}$, $\ m=1,2$, where
\[
\lb{idm} \m_{m0}=\cos mz,\ \ \ \m_{m1}={\sin mz\/4z}mV_{(1)},\ \ \
\m_{m2}={1\/8z^2} (I_m^0(\l)-{m^2\/2}\cos mz V_{(2)}),
\]
and $z=\sqrt \l$, where $V_{(m)}$ is defined by \er{dvm}. We
obtain $\m_1^2=B_1+B_2$, where
$$
 B_1=(\m_{10}+\m_{11})^2+2\m_{10}\m_{12},\ \
B_2=2\m_{10}\wt\m_{13}+
(\m_1-\m_{10}+\m_{11})(\m_1-\m_{10}-\m_{11}).
$$
Then \er{idm} yields
$$
B_1=\cos^2 z+{\sin 2z\/4z}V_{(1)}+{4\cos zI_1^0(\l)+\sin^2
zV_{(1)}^2-2\cos^2 z V_{(2)}\/16z^2}.
$$
Thus we get
$$
\r={\m_2+1\/2}-\m_1^2=G_1+G_2,\ \ \
G_1={1+\m_{20}+\m_{21}+\m_{22}\/2}-B_1,  \ \ \ G_2={\wt
\m_{23}\/2}-B_2,
$$$$
G_1={I_2^0(\l)-2\cos 2z V_{(2)}-4\cos zI_1^0(\l)+2\cos^2 z
V_{(2)}-\sin^2
zV_{(1)}^2\/16z^2}={\sin^2z\/16z^2}(2V_{(2)}-V_{(1)}^2),
$$
which yields $G_1=\r_0^2$. Using \er{em01}-\er{em}  we obtain
$$
|B_2|\le 2{\vk^3\/3!}A^2+(\vk A+\vk
A){\vk^2\/2!}A^2={4\vk^3\/3}A^2,
$$
\[
|G_2|\le {|\wt \m_{m3}|\/2}+|B_2|\le {4(\vk)^3\/
3!}A^2+{4\vk^3\/3}A^2=2\vk^3A^2,
\]
which yields \er{3a}.

 \no ii) Let $N'>N$ be another integer and $r=\pi
(N+{1\/2})$. Note that
$\vk\le{V_0^2\/2^7\pi\|V\|_1^2}\le{1\/2^7\pi}$ for $\l \in \cD_r$.
Using the estimate $e^{|\Im \sqrt\l|}<4|\sin \sqrt\l|,\l\in \cD_r$
and \er{3a} we obtain (on all contours)
\[
\lb{er1} |\r(\l)-\r_0^2(\l)|\!\le\! 2{\|V\|_1^3\/|z|^3}e^{2(|\Im
z|+\vk)}\!\le\! {C_0\/\sqrt{|\l|}}|\r_0^2(\l)|,\ \ \ \ C_0\ev
4^4{\|V\|_1^3\/c_0^2}<{\sqrt{|\l|}\/2},\ \l\in \cD_r.
\]
Hence, by the Rouch\'e theorem, $\r(\l)$ has as many roots,
counted with multiplicity, as ${\sin^2 \sqrt\l\/\l}$ in the
bounded domain $\cD_r\cap \{\l:|\l|\le \pi (N'+{1\/2}) \}$.
Since ${\sin^2 \sqrt\l\/\l}$ has exactly one double root at $(\pi
n)^2,n\ge 1$, and since $N'>N$ can be chosen arbitrarily large,
 ii) follows.

\no iii) Let $\r\ev\r_0^2+\r_1,\ \r_0=c_0{\sin
\sqrt\l\/2\sqrt\l}$. Estimates \er{er1}  imply
\[
\lb{3ae}
 \sqrt{\r(\l)}=\r_0(\l)\sqrt{1+b(\l)}, \ \ \
b={\r_1\/\r_0^2},\ \ \ |b(\l)|\le {C_0\/|z|}<{1\/2},\ \ \l\in
\cD_r.
\]
Using the estimate $|\sqrt{1+b(\l)}-1|\le {3\/5}|b(\l)|$ for
$\l\in \cD_r$, we obtain
$$
|\sqrt{\r(\l)}-\r_0(\l)|=|\r_0(\l)(\sqrt{1+b(\l)}-1)|\le
{3C_0\/5}{|\r_0(\l)|\/\sqrt{|\l|}}.\ \ \ \ \BBox
$$

We need asymptotics of the eigenvalues of the monodromy
matrix.

\no {\bf Lemma 3.2.} {\it Let $V\in \cH^1, V_0>0$. Then the
monodromy matrix $M(\l)$ has two multipliers $\t_m(\l),m=1,2$ with
asymptotics \er{1ca}. Moreover, asymptotics \er{asL} is fulfilled.

 \no Proof.} Using estimates \er{em2}, \er{3ad} we obtain
asymptotics \er{asL}. Asymptotics \er{asL} yield
$$
\D_m^2(\l)-1=\lt(\cos \sqrt\l+V_{m0}{\sin
\sqrt\l\/2\sqrt\l}\rt)^2-1+O(E(\l))
$$
$$
=-\sin^2 \sqrt\l+V_{m0}{\sin \sqrt\l\cos \sqrt\l\/\sqrt\l}+
O(E(\l))=-\sin^2 \sqrt\l \lt(1-V_{m0}{\cos \sqrt\l\/\sqrt\l
\sin
\sqrt\l}+O(E(\l))\rt)
$$
as $\l\in\cD_1, |\l|\to \iy$, where $E(\l)={e^{2|\Im \sqrt
\l|}\/\l}$, which implies
\[
\sqrt{\D_m^2(\l)-1}=i\sin \sqrt\l -iV_{m0}{\cos
\sqrt\l\/2\sqrt\l}+O(E(\l)).
\]
By \er{2l}, the matrix $M(\l),\l\in\cD_r$, has the eigenvalues
$\t_m(\l)$ satisfying the identities
$\t_m(\l)+\t_m(\l)^{-1}=2\D_m(\l)$. Then $\t_m(\l)$ has the form $
\t_m(\l)=\D_m(\l)+\sqrt{\D_m^2(\l)-1} $ and the asymptotics give
$$
\t_m(\l)=e^{i\sqrt\l}+V_{m0}{\sin \sqrt\l\/2\sqrt\l} -iV_{m0}{\cos
\sqrt\l\/\sqrt\l}+O(E(\l))=e^{i\sqrt\l}\lt(1-{iV_{m,0}\/2\sqrt
\l}\lt)+O(E(\l))
$$
which yields \er{1ca}.
 $\BBox$

Now we prove our first result about the Lyapunov function 
$\D=\m_1+\sqrt \r$.

\no{\bf Proof of Theorem 1.1.} From Lemma 3.1 we obtain the
analytic function $\D$ on the Riemann surface of the function
$\sqrt \r$.

\no  i) Identity \er{2l} shows that $\D_m={\t+\t^{-1}\/2}$ for
some multiplier $\t$. Lemma 3.2 gives the asymptotics of $\D_m$
and $\t_m, m=1,2$.

\no ii) The result follows from the statement i) and the Lyapunov
Theorem (see Sect.1).

\no iii) If $\l\in\s(T)$, then $\m_1(\l)$ is real. By ii), $\D(\l)$ is
also real. Hence by \er{1g}, $\sqrt{\r(\l)}$ is real and $\r(\l)\ge 0$.

\no iv) Assume that $\D_m'(\l_0)=0$ for some $\l_0\in I\ss \s(T)$
and $m\in \{1,2\}$. Then we have
\[
\lb{6a}
\D_m(\l)=\D_m(\l_0)+(\l-\l_0)^p{\D_m^{(p)}(\l_0)\/p!}+O(|\l-\l_0|^{p+1}),\
\ \ {\rm as} \ \ \ \l-\l_0\to 0,
\]
where $\D_m^{(p)}(\l_0)\ne 0$ for some $p\ge 1$. By the
Implicit Function Theorem, there exists some curve $\G\ss
\{\l:|\l-\l_0|<\ve \}\cap \C_+, \G\neq \es$,  for some $\ve>0$
such that $\D_m(\l)\in (-1,1)$ for any $\l\in \G$. Thus we have
a contradiction with the Lyapunov Theorem in Sect. 1.
 $\BBox$

Recall that $D_{\pm}(\l)={1\/4}\det(M(\l)\mp I_4)$, the set
$\{\l:D_+(\l)=0\}$ is a periodic spectrum and the set
$\{\l:D_-(\l)=0\}$ is an anti-periodic spectrum. Now we prove a
lemma about the number of periodic and anti-periodic eigenvalues in
a large disc.

\no {\bf Lemma 3.3.} {\it For each $V\in \cH^1$ the functions
$\D_1+\D_2,\D_1\D_2, D_\pm$ are entire and satisfy the following
identities:
\[
\lb{3f}
 \D_1^2+\D_2^2=1+\m_2, \ \ \ \D_1\D_2=2\m_1^2-{\m_2+1\/2},
\]
\[
\lb{5d}
 D_{\pm}=(\m_1\mp 1)^2-\r={(2\m_1\mp1)^2-\m_2\/2},\ \
\ D_+-D_-=-4\m_1.
\]
Let in addition $\Tr \int_0^1V(t)dt=0$. Then the following
estimates and properties are fulfilled:
\[
\lb{3i}
 \max\lt\{\lt|D_+(\l)-4\sin^4{\sqrt\l\/2}\rt|,
 \ \ \ \lt|D_-(\l)-4\cos^4{\sqrt\l\/2}\rt|\rt\}
\le{\vk^2}(2+\vk)^2e^{2|\Im\sqrt\l|+2\vk}, \ \ \ \l\in\C.
\]
\no i) For each integer $N>8\|V\|_1$ the function $D_+$ has
exactly $4N+2$ roots, counted with multiplicity, in the disc
$\{|\l|<4\pi^2(N+{1\/2})^2\}$ and for each $n>N$, exactly four
roots, counted with multiplicity, in the domain $\{|\sqrt\l-2\pi
n|<{\pi\/2}\}$. There are no other roots.

\no ii) For each integer $N>8\|V\|_1$ the function $D_-$ has
exactly $4N$ roots, counted with multiplicity, in the disc
$\{|\l|<4\pi^2 N^2\}$ and for each $n>N$, exactly four roots,
counted with multiplicity, in the domain $\{|\sqrt\l-\pi
(2n+1)|<{\pi\/2}\}$. There are no other roots.

\no Proof.} By Lemma 2.1, 2.2 the functions $\D_1+\D_2,\D_1\D_2,
D_\pm$ are entire and identities \er{3f}, \er{5d} are fulfilled.
Using Lemmas 2.1 and 3.1 we obtain
$$
\m_1(\l)=\cos z+\wt\m_{12}(\l),\ \ |\wt\m_{12}(\l)|\le
{\vk^2\/2}A^2,\ \ |\r(\l)|\le (1+2\vk){\vk^2\/2}A^2.
$$
Substituting these relations into $D_{\pm}=(\m_1\mp 1)^2-\r$ we
get \er{3i}.

\no i) Let $N'>N$ be another integer. Let $\l$ belong to the
contours $C_0(2N+1),C_0(2N'+1),C_{2n}({1\/2}),|n|>N$, where
$C_n(r)=\{\sqrt\l:|\sqrt\l-\pi n|=\pi r\},r>0$.  Note that $\vk\le
{1\/16\pi}$ on all contours. Then \er{3i} and the estimate
$e^{{1\/2}\Im \sqrt\l}<4|\sin{\sqrt\l\/2}|$ on all contours yield
$$
\lt|D_+(\l)-4\sin^4{\sqrt\l\/2}\rt|\le {e^{2\vk}(2+\vk)^2\/4^4
\pi^2}e^{2|\Im \sqrt\l|}<{1\/4}\lt|4\sin^4{\sqrt\l\/2}\rt|.
$$
Hence, by Rouch\'e's theorem, $D_+(\l)$ has as many roots, counted
with multiplicities, as $\sin^4{\sqrt\l\/2}$ in each of the
bounded domains and the remaining unbounded domain. Since
$\sin^4{\sqrt\l\/2}$ has exactly one root of the multiplicity four
at $(2\pi n)^2$, and since $N'>N$ can be chosen arbitrarily large,
i) follows. The proof for $D_-$ is similar.
 $\BBox$

We are ready to prove Theorem 1.2.

\no {\bf Proof of Theorem 1.2.} i) It is enough to consider the
case $V_{10}=c_0=-V_{20}$.  We prove asymptotics \er{a8} for
$\l_{2n,m},1\le m\le 4$. The proof for $\l_{2n+1,m}$ is similar.
Firstly, we prove the rough asymptotics of $\l_{2n,m},r_n$. Lemma
3.3.i yields $\sqrt{\l_{2n,m}}=2\pi n+\ve_n, |\ve_n|<1$, as $n\to
\iy$. By Lemma 3.3, $D_+(\l)=4\sin^4{\sqrt\l\/2}+O(n^{-2})$ as $
|\sqrt\l-\pi n|\le 1,n\to\iy$. Then identity $D_+(\l_{2n,m})=0$
implies
\[
\lb{3p}
\sqrt{\l_{2n,m}}=2\pi n+\ve_n,\ \ \
\ve_n=O(n^{-1/2}),\ \ \ 1\le m\le 4.
\]
Lemma 3.1.ii implies $\sqrt{r_{2n-m}}=\pi n+\d_n,|\d_n|<1$ for
$n\to \iy,m=0,1$. Moreover, Lemma 3.1.i gives
$\r(\l)=c_0^2{\sin^2\sqrt\l\/4\l}+O(n^{-3}), |\sqrt\l-\pi n|\le 1$
as $n\to\iy$. Since $\r(r_n)=0$, we have
\[
\lb{3r}
 \sqrt{r_{2n-m}}=\pi n+\d_n,\ \ \ \d_n=O(n^{-1/2}),\ \ \ m=0,1.
\]
Secondly, in order to improve these asymptotics of $\l_{2n,m},r_n$
we need asymptotics of the multipliers in a neighborhood of the
points $\pi n$. We introduce the matrix $\wt M=U^{-1}MU$ with the
same eigenvalues, where $U=\ma I_2& 0\\ 0& \sqrt\l I_2\am$. We
shall show the asymptotics
\[
\wt M(\l)=\wt M_0(\l)+\ell_1^2(n),\ \ \wt M_0(\l)=\ma C(\l)& S(\l)\\ -S(\l)&
C(\l)\am, \ \ \ \ \sqrt\l=\pi n+O(n^{-1/2}),\ \ \ \lb{3l}
\]
where
\[
\lb{3k}
C(\l)=\diag(\cos\e_+,\cos\e_-),\ \ \
S(\l)=\diag({\sin\e_+},{\sin\e_-}),\ \ \ \e_{\pm}=\sqrt{\l\mp
c_0}. 
\]
Estimate \er{2if} gives
$$
\vt(1,\l)=\cos\sqrt\l I_2+{1\/\sqrt\l}\int_0^1\sin\sqrt\l(1-t)\cos\sqrt\l tV(t)dt+
O(n^{-2})
$$
\[
\lb{2e} =\cos\sqrt\l
I_2+{1\/2\sqrt\l}\int_0^1\lt(\sin\sqrt\l+\sin\sqrt\l(1-2t)\rt)V(t)dt+O(n^{-2}),
\ \ \ |\sqrt\l-\pi n|\le 1,
\]
as $n\to+\iy$. Let $\sqrt\l=\pi n+u_n,u_n=O(n^{-1/2})$.  The
Taylor formula gives $\sin 2t(\pi n+u_n)=\sin 2\pi nt+2tu_n\cos
2\pi nt+O(n^{-1})$ and the similar formula for $\cos 2t(\pi
n+u_n)$. Substituting these asymptotics into \er{2e}  we obtain
$$
\vt(1,\l)=\cos\sqrt\l I_2+{\sin\sqrt\l\/2\sqrt\l}c_0J+\ell_1^2(n)=C(\l)+\ell_1^2(n),
\ \ \ \sqrt\l=\pi n+O(n^{-1/2}).
$$
Similar arguments for $\vp(1,\l),\vt'(1,\l),\vp'(1,\l)$ yield
$$
\vp(1,\l)={1\/\sqrt\l}S(\l)+\ell_2^2(n),\ \ \ \vt'(1,\l)=-\sqrt\l
S(\l)+\ell^2(n),\ \ \ \vp'(1,\l)=C(\l)+\ell_1^2(n),
$$
as $\sqrt\l=\pi n+O(n^{-1/2})$. Substituting the obtained
asymptotics into the definition \er{1c} of $M$ and using the identity $\wt M=U^{-1}MU$
we get \er{3l}.

We will use standard arguments from the perturbation theory
(see [Ka,p.291]). Let $A,B$ be bounded operators, $A$ be a normal
operator and $\s(A),\s(B)$ be spectra of $A,B$. Then
$\dist\{\s(A),\s(B)\}\le\|A-B\|$. Note that $\wt M_0(\l)$ is a normal
operator having eigenvalues $e^{i\e_{\pm}}$. Hence $\wt M$ has
eigenvalues $\t_{\pm}$ satisfying the estimates
$|\t_{\pm}-e^{i\e_{\pm}}|<|\wt M-\wt M_0|$. Then \er{3l} implies
\[
\t_{\pm}(\l)=e^{i\e_{\pm}}+\ell_1^2(n),\ \ \ {\rm as}\ \ \ \sqrt\l=\pi
n+O(n^{-1/2}). \lb{3m}
\]

Now we improve asymptotics \er{3p} for $\l_{2n,m}$. Note that
$\l=\l_{2n,m}$  iff $\t_+(\l)=1$ or $\t_-(\l)=1$. Then  \er{3m}
and \er{3p} yield
\[
\lb{3n} e^{i\e_{\pm}(\l_{2n,m})}=1+\ell_1^2(n),\ \ \ \l=\l_{2n,m}.
\]
Substituting  asymptotics \er{3p} into \er{3n} we have
$$
\sqrt{(2\pi n+\ve_n)^2\pm c_0}=2\pi n+\ell_1^2(n),\ \ \ n\to+\iy.
$$
and therefore $\ve_n=\mp{c_0\/4\pi n}+\ell_1^2(n)$. Substituting this
asymptotics into \er{3p} we obtain \er{a8} for $\l_{2n,m}$.

We improve asymptotics \er{3r} of $r_n$. Note that $\l=r_{n}$, iff
the following condition is fulfilled: $\t_+(\l)=\t_-^{-1}(\l)$ or
$\t_+(\l)=\t_-(\l)$. Asymptotics \er{3m} imply that the second
equation has no solutions for large $n$. We rewrite the first
equation in the form:
\[
 \lb{3o}
e^{i\e_+}=e^{-i\e_-}+\ell_1^2(n),\ \ \ {\rm as}\ \ \ \sqrt\l=\pi
n+O(n^{-1/2}).
\]
Substituting asymptotics \er{3r} into \er{3o} we obtain
$$
\sqrt{(\pi n+\d_n)^2-c_0}+\sqrt{(\pi n+\d_n)^2+c_0}=2\pi
n+\ell_1^2(n),
$$
which gives $\d_n=\ell_1^2(n)$. Substituting this asymptotics into \er{3r} we obtain
\er{a9}.

\no ii) a) Assume that we have the
periodic spectrum $\l_{0,m},m=1,2,\l_{2n,m},m=1,2,3,4,n\ge 1$.
Using the asymptotics \er{a8} and repeating the standard
arguments (see [PT,pp.39-40]) we obtain the Hadamard factorization
for the function $D_+$:
$$
D_+(\l)={(\l-\l_{0,1})(\l-\l_{0,2})\/4} \prod_{m=1,2,3,4,n\ge
1}{\l_{2n,m}-\l\/(2\pi n)^2}.
$$
In a similar way, we determine $D_-$ by the anti-periodic
spectrum. Using \er{5d} we obtain $\r$. Thus, we recover the
resonances.

 b) Suppose, that we have the periodic spectrum
and the set of the resonances.  Then we determine the functions
$D_+$ by the periodic spectrum and $\r$ by the resonances. Using
\er{5d} we get $\m_1$, $\m_2$ and then $D_-$. Thus, we recover the
anti-periodic spectrum. The proof of another case is similar.
$\BBox$

\section {Example}
\setcounter{equation}{0}

\no {\bf 1. Periodic $\d$-potentials.} Consider the operator
$T^\g=-{d^2\/dx^2}+q^\g, \g\in \C$, where $q^\g=a J+\g \d_{per} J_1$ and the potential $\d_{per}=\sum_{-\iy}^\iy \d(x-n-{1\/2})$. Let $\m_1^\g=\m_1(\cdot,q^\g),\r^\g=\r(\cdot,q^\g),...$

\no {\bf Lemma 4.1.} {\it For the operator $ T^\g=-{d^2\/dx^2}+a
J+\g \d_{per} J_1$ the following identities are fulfilled:
\[
\lb{idg} \m_1^\g=\m_1^0,\ \ \ \m_2^\g=\m_2^0+2h,\ \ \
\r^\g=\r^0+h,\ \ \ D_{\pm}^\g=D_{\pm}^0-h,\ \ \
h={\g^2\/4}s_+s_-,\ \ \ \ \ \g\in \C.
\]
\no Proof.} The solution $y(x), x\in \R$ of the system $-y''+q^\g
y=\l y$ is continuous and $y'(x_n+0)-y'(x_n-0)=\g J_1 y(x_n)$  at
the points $x=x_n=n+{1\/2}$. Then the fundamental solutions have
the form: $ \vt^\g(x)=\vt^0(x),\ \  \vp^\g(x)=\vp^0(x),  \ \ 0\le
x <{1\/2},
$ and
$$
\vt^\g(x)=\ma c_+(x)&\g c_-({1\/2})s_+(x-{1\/2})\\
\g c_+({1\/2})s_-(x-{1\/2})&c_-(x)\\ \am, \ \ \ {1\/2}\le x
<{3\/2},
$$
$$
\vp^\g (x)=\ma s_+(x)&\g s_-({1\/2})s_+(x-{1\/2})\\
\g s_+({1\/2})s_-(x-{1\/2})&s_-(x)\\ \am ,\ \ \ {1\/2}\le x
<{3\/2},
$$
and
$$
\vt^\g(x)=\ma  c_+(x)+\g^2c_+({1\/2})s_-(1)s_+(x-{3\/2}) & *\\
 * & c_-(x)+\g^2c_-({1\/2})s_+(1)s_-(x-{3\/2})\\ \am,\
$$$$
\vp^\g (x)=\ma s_+(x)+\g^2s_+({1\/2})s_-(1)s_+(x-{3\/2}) & *\\
 * & s_-(x)+\g^2s_-({1\/2})s_+(1)s_-(x-{3\/2})\\
\am,
$$
for ${3\/2}\le x <{5\/2}$, where  $*$ is some term.  These
relations yield
$$
\vt(1)=\ma c_+(1)&*\\ *&c_-(1)\\ \am,\ \ \ \vp'(1)=\ma c_+(1)&*\\
*&c_-(1)\\ \am,
$$
$$
\vt(2)=\ma c_+(2)+2h& *\\  * & c_-(2)+2h\\ \am, \ \ \ \vp'(2)=\ma
 c_+(2)+2h & *\\ * & c_-(2)+2h\\ \am.
$$
The last identities and \er{m0} imply
$$
\m_1^\g={\Tr (\vt^\g(1)+(\vp^\g)'(1))\/4}=\m_1^0,\ \ \ \ \ \ \
\m_2^\g={\Tr (\vt^\g(2)+(\vp^\g)'(2))\/4} =\m_2^0+2h,
$$
which give $\r^\g={1\/2}(1+\m_2^\g)-(\m_1^\g)^2=\r^0+h$ and $
D_{\pm}^\g=(\m_1^\g\mp 1)^2-\m_2^\g=D_{\pm}^0-h$. $\BBox$

We describe the spectrum of the operator $T^\g$.

\no {\bf Lemma 4.2.} {\it Let the operator $
T^\g=-{d^2\/dx^2}+q^\g$, where $q^\g=a J+\g \d_{per} J_1, a>0$ and $\g\in\R$.

\no i) Let $0<{a\/2\pi^2}-n_a<1$, for some integer $n_a\ge 0$.
Then for each $n\in \N$ there exist analytic functions
$z_n^{\pm}(\g), |\g|<\g_n$ for some $\g_n>0$  such that
$z_n^{\pm}(\g)$ is the zero of $\r^\g(\l)$ and
\[
z_n^{\pm}(\g)=z_n^0\pm (-1)^n\g i^{k_n}(c_n+O(\g)),\ \  c_n>0,\ \
k_n=\cases {1\ \ &if\ \ $1\le n\le n_a,$\cr
          0\ \ &if\ \  $n>n_a$\cr}, \ \ \g\to 0.
\]
Moreover, each spectral interval $(z_n^-(\g),z_n^+(\g))\ss\R,n>n_a$ is a gap in the spectrum of $T^\g$.

 \no ii) If ${2a\/\pi^2}\notin \N$, then for each $n\ge
0, m=1,2$ there exist real analytic functions $\l_{m,n}^{\pm}(\g),
\g\in (-\g_n,\g_n)$ for some $\g_n>0$  such that
$\l_{n,m}^{\pm}(\g)$ is the zero of the function $\D_{m}^2(\l)-1$
and
\[ \l_{n,m}^{-}(\g)<\l_{n,m}^{+}(\g),\ \ \ \ \ \ \
\l_{n,m}^{\pm}(0)=\l_{n,m}^0.
\]
Moreover, each spectral interval
$(\l_{n,m}^{-}(\g),\l_{n,m}^{+}(\g))\neq \es, n\ge 1$ has
multiplicity 2.

 \no Proof.} i) Recall that $\e_{\pm}=\sqrt{\l\mp a}$. The zero of $\r^\g(\l)=0$ satisfies the equation
 \[
0=\r^\g(\l)=f^2(\l)+{\g^2\/4}s_+(\l)s_-(\l),\ \ \ \ \ f\ev
2\sin{\e_+-\e_-\/2}\sin{\e_++\e_-\/2}.
 \]
The zeros $z_n^0$ of $\r^0=\r(\cdot, q^0)$ have  the form \er{4c}
and satisfy the following identities
 \[
\sqrt{z_n^0+a}+\sqrt{z_n^0-a}={a\/\pi n},\ \ \ \ \ \
\sqrt{z_n^0+a}-\sqrt{z_n^0-a}=2\pi n,\ \ {\rm if}\ \ 1\le n\le
n_a, \lb{4d}
\]
\[
\lb{4e} \sqrt{z_n^0+a}+\sqrt{z_n^0-a}=2\pi n,\ \ \ \ \ \ \ \
\sqrt{z_n^0+a}-\sqrt{z_n^0-a}={a\/\pi n}, \ \ \ \ \ \ {\rm if}\ \
n> n_a\ge 0.
\]
Using \er{4d}, \er{4e} we have the identity
\[
\sin \e_+(\l)\sin \e_-(\l)|_{\l=z_n^0}=\cases {1-\cos {a\/\pi
n}>0\ \ \ &if\ \ \ \ $1\le n\le n_a,\ \ \ \ $\cr
          \cos {a\/\pi n}-1<0\ \ \ &if\ \ \ \ $n>n_a$\cr}.
\]
Recall that the function $f$ has only simple zeros $\l=z_n^0, n\ge
1$.
 Consider the case $n>n_a$, the proof
for $1\le n\le n_a$ is similar. We rewrite $\r^\g(\l)=0$ in the
form
$$
(f(\l)-\g F(\l))(f(\l)+\g F(\l))=0,\  \ \ \ \ \ F(\l)\ev\sqrt
{-s_+(\l)s_-(\l)/2},\ \ F(z_n^0)>0.
$$
Applying the Implicit Function Theorem to $\F_{\pm}(\l,\g)=0$,
where $\F_{\pm}(\l,\g)=f(\l)\pm\g F(\l)$ and ${\pa\/\pa
\l}\F_{\pm}(z_n^0,0)\neq 0$ we obtain a unique solution
$z_n^{\g,\pm}=u=u(\g), \g\in (-\g_0,\g_0), u(0)=z_n^0$ of the
equation $\F(\l,\g)=0$, such that $\F(u(\g),\g)=0,\g\in
(-\g_0,\g_0)$ for some $\g_0>0$.

ii) Consider the equation
$0=D_+^\g(\l)=D_+^0(\l)-{\g^2\/4}s_+(\l)s_-(\l)$.
Using
\[
D_+^0(\l)=4\sin^2{\e_+\/2}\sin^2{\e_-\/2},\ \ \ \ \ \ \
s_{\pm}(\l)={2\/\e_{\pm}}\sin{\e_{\pm}\/2}\cos{\e_{\pm}\/2},\ \
\cos \e_-(\l_{2n,1})\neq 1
\]
we obtain 2 equations for the zeros of $D_+^\g(\l)$
 \[
 \lb{eq4}
\F(\l,\g)\ev\sin{\e_+\/2}-{\g^2\/4\e_-\e_+}\cos{\e_+\/2}\cot
{\e_-\/2}=0,\ \ \ \ \  \ \ \ \sin{\e_+\/2}=0.
 \]
 The equation ${\sin{\e_+\/2}\/\e_+}=0$ has the zeroes
 $\l_{1,2n}^0$.
Consider the first equation in \er{eq4}.
 Applying the Implicit Function Theorem to $\F(\l,\g)=0$, where
 ${\pa\/\pa \l}\F(\l_{2n,1}^0,0)\neq 0$ we obtain a unique
solution $u=u_n(\g), \g\in (-\g_n,\g_n), u_n(0)=\l_{2n,1}^0$ of
the equation $\F(\l,\g)=0$, such that $\F(u_n(\g),\g)=0,\g\in
(-\g_n,\g_n)$ for some $\g_n>0$. The proof for $D_-^\g(\l)=0$ is
similar.$\BBox$

\no {\bf 2. The perturbed operator.}
 We consider the operator $T_{\g,\n}=-{d^2\/dx^2}+q_{\g,\n}$ where the potential $q_{\g,\n}=aJ+\g v_\n J_1$  satisfies the Condition A, $\g\in \R$ is small and $a>0$. We determine the asymptotics of the function $\r(\l,q_{\g,\n}), \m_m(\l,q_{\g,\n}), m=1,2$.

\no {\bf Lemma 4.3.} {\it Each function
$\r(\l,q_{\g,\n}),\m_m(\l,q_{\g,\n}),m=1,2,\n=1,{1\/2},{1\/3},...,$
is analytic in $\C^2$. Moreover, uniformly on any compact in
$\C^2$ the following asymptotics are fulfilled:
\[
\lb{4a}
  \r(\l,q_{\g,\n})=\r^\g(\l)+o(1), \ \ \ \ \
 \m_m(\l,q_{\g,\n})=\m_m^\g(\l)+o(1), m=1,2,\ \ \  \n\to 0.
\]
\no Proof.} The fundamental solutions $\vt_{\n,\g},\vp_{\n,\g}$ of
the equation $-y''+q_{\g,\n}y=\l y$, satisfying the conditions
$\vt_{\n,\g}(0,\l)=(\vp_{\n,\g})'(0,\l)=I_2$,
$(\vt_{\n,\g})'(0,\l)=\vp_{\n,\g}(0,\l)=0$ are the solutions of
the integral equation
\[
\lb{fd} \vp_{\n,\g}(x,\l)\!=\vp_0(x,\l)+\int_0^x\!\!\vp_0(x-t,\l)
q_{\g,\n}(t)\vp_{\n,\g}(t,\l)dt.
\]
The standard iterations in \er{fd} yield
\[
\lb{4id} \vp_{\n,\g}={\sum}_{n\ge 0} \vp_{n,\n,\g}\,, \quad
\vp_{n,\n,\g}(x,\l)= \int_0^x\vp_0(x-t,\l)
q_{\g,\n}(t)\vp_{n-1,\n,\g}(t,\l)dt.
\]
The last identity gives
\[
\vp_{n,\n,\g}(x,\l)= \int\limits_{0< x_1< x_2< ...<
x_{n+1}=x}\!\!\!\!\vp_0(x_1,\l) \lt(\prod\limits_{1\le k\le
n}^{\curvearrowleft} \vp_0(x_{k+1}-x_k,\l)
q_{\g,\n}(x_k)\rt)dx_1dx_2...dx_n, \lb{4ig}
\]
where for matrices $a_1,a_2,...,a_n$ we define $\prod\limits_{1\le
k\le n}^{\curvearrowleft}a_k=a_na_{n-1}...a_1$. Substituting the
estimate $\|\sqrt\l\vp_0(x,\l)\|\le e^{|\Im \sqrt\l|x}$ into
\er{4ig} we obtain $\| \sqrt\l \vp_{n,\n,\g}(x,\l)\|\le
{(2x(a+|\g|))^n\/n!}e^{|\Im z|x}, $ which shows that for each
$x\ge 0$ the formal series \er{4id} converges uniformly on bounded
subsets of $\C$. Each term of this series is an entire function.
Hence the sum is an entire function.
Since $v_\n\to p$ in the sense of distributions we obtain
$\vp_{n,\n,\g}(x,\l)\to \vp_{n,0,\g}(x,\l)$ as $\n \to 0$
uniformly on bounded subset of $\R\ts\C^2$, which yields \er{4a}.
 $\BBox$

\no We give

\no {\bf Proof of Proposition 1.3.} Lemma 4.3 yields
$\r(\l,q_{\g,\n})\to \r(\l,q^\g,)$ and $\m_m(\l,q_{\g,\n})\to
\m_m(\l,q^\g),m=1,2,$  uniformly on any compact set in $\C^2$ as
$\n\to 0$. Then their zeros converge to the corresponding zeros
at $\n=0$, uniformly on any compact in $\C^2$. Due to Lemma 3.3 we
have convergence of $D_\pm (\l,q_{\g,\n})$ and of the Lyapunov
function $\D(\l,q_{\g,\n})$. Thus using Lemmas 4.1-4.2 we obtain
the proof of Proposition 1.3.$\BBox$

\no {\bf References}

\no [BKW] Brown, B., Knowles, I., Weikard, R.: On the inverse resonance problem. J. London Math. Soc. (2) 68 (2003), no. 2, 383-401

\no [Ca1] Carlson, R.: Eigenvalue estimates and trace formulas for
the matrix Hill's equation. J. Differential Equations 167 (2000),
no. 1, 211--244.

\no [Ca2]  R. Carlson. Large eigenvalues  and trace formulas for
the matrix Sturm-Liouville problems. SIAM J. Math. Anal.30 (1999),
949-962

\no [Ca3] Carlson, R.: Compactness of Floquet isospectral sets for
the matrix Hill's equation. Proc. Amer. Math. Soc. 128 (2000),
 no. 10, 2933--2941.

\no [CHGL] Clark S., Holden H., Gesztesy, F., Levitan, B.:
Borg-type theorem for matrix-valued Schr\"odinger and Dirac
operators, J. Diff. Eqs. 167(2000), 181-210

\no [DS] Dunford, N. and Schwartz, J.: Linear Operators Part II:
Spectral Theory, Interscience, New York, 1988.

\no [GL] Gel'fand I., Lidskii, V.: On the structure of the regions
of stability of linear canonical
   systems of differential equations with periodic coefficients. (Russian)
   Uspehi Mat. Nauk (N.S.) 10 (1955), no. 1(63), 3--40.

\no [Ge] Gel'fand, I.:
 Expansion in characteristic functions of an equation with periodic coefficients.
   (Russian) Doklady Akad. Nauk SSSR (N.S.) 73, (1950). 1117--1120.

\no [KK] Kargaev P., Korotyaev E.: Inverse problem for the Hill
operator, the direct approach.  Invent. Math., 129(1997), no. 3,
567-593

\no [Ka] Kato, T.: Perturbation theory for linear operators.
Springer-Verlag, Berlin, 1995.

\no [K1] Korotyaev, E.: Characterization of the spectrum of
Schr\"odinger operators with periodic distributions. Int. Math.
Res. Not.  (2003) no. 37, 2019--2031

\no [K2] Korotyaev, E.: The estimates of periodic potentials in
terms of effective masses. Comm. Math. Phys.
   183 (1997), no. 2, 383--400.

\no [K3] Korotyaev, E.:  Stability for inverse resonance problem, Int. Math. Res. Not. 73(2004), 3927-3936

\no [Kr] Krein, M.: The basic propositions of the theory of
$\lambda$-zones of stability of a canonical system of linear
differential equations with periodic coefficients. In memory of A.
A. Andronov, pp. 413--498. Izdat. Akad. Nauk SSSR, Moscow, 1955.

\no [Ly] Lyapunov, A.: The general problem of stability of motion,
2 nd ed. Gl. Red. Obschetekh. Lit., Leningrad, Moscow, 1935;
reprint Ann. Math. Studies, no. 17, Prinston Univ. Press, Prinston, N.J. 1947

\no [PT] P\"oshel,J., Trubowitz, E.: Inverse spectral theory. Pure
and Applied Mathematics, 130. Academic Press, Inc., Boston, MA,
1987. 192 pp.

\no [RS] Reed, M., Simon, B.: Methods of Modern Mathematical
Physics, Vol.II, Fourier Analysis, Self-Adjointness, Academic Press, New York, 1975

\no [S] Simon B.: Resonances in one dimension and Fredholm determinants,
 J. Funct. Anal. 178 (2000), no. 2, 396-420.

\no [YS] Yakubovich, V., Starzhinskii, V.: Linear differential
equations with periodic coefficients. 1, 2.
   Halsted Press [John Wiley \& Sons] New York-Toronto,
   1975. Vol. 1, Vol. 2

\no [Z] Zworski M.: Distribution of poles for scattering
on the real line, J. Funct. Anal. 73, 277-296, 1987

\end{document}